\documentclass[12pt,a4paper]{article}
 
\usepackage{graphicx,color,url}
\usepackage{xspace}
\usepackage{caption}

\usepackage{amsfonts,amsmath,amssymb,amsthm}
\usepackage{latexsym}
\usepackage{ulem} 
\usepackage[T1]{fontenc}
\usepackage{enumitem}
\usepackage[toc,page]{appendix}

\newtheorem{theorem}{Theorem}[section]
\newtheorem{definition}[theorem]{Definition}
\newtheorem{example}[theorem]{Example}
\newtheorem{lemma}[theorem]{Lemma}
\newtheorem{proposition}[theorem]{Proposition}

\newtheorem{notation}[theorem]{Notation}

\title{\Large Domination, matching and transversal numbers for Berge-$G$  hypergraphs\footnote{Partially supported by PAI-DI-FQM 189, PAI-DI-FQM 326, (Andalucía Goverment, Spain) and VII Plan de Investigación y Transferencia,  SOL2024-30793,    Universidad de Sevilla, Spain.}}

\author{
\large Mar\'{\i}a Jos\'e Ch\'avez de Diego\footnote{Dpto. de Matemática Aplicada I, Universidad de Sevilla, Spain, mjchavez@us.es.},\\
Pablo Montero Moreno\footnote{Dpto. de Geometr\'{\i}a y  Topolog\'{\i}a, Universidad de Sevilla, Spain, pabmonmor1@alum.us.es}\\
and Mar\'{\i}a Trinidad Villar-Liñ\'an\footnote{Dpto. de Geometr\'{\i}a y  Topolog\'{\i}a, Universidad de Sevilla, Spain, villar@us.es } }

\begin{document}
\maketitle

\begin{abstract}
Let $G=(V(G),E(G))$ be a graph and $H=(V(H),E(H))$ be a hypergraph. The hypergraph $H$ is a {\it Berge-G} if there is a bijection $f : E(G) \mapsto E(H)$  such that for each $e \in  E(G)$  we have $e \subseteq  f(e)$.  We define {\it dilations of $G$} as a particular subfamily of not necessarily uniform Berge-$G$ hypergraphs. We examine domination, matching and transversal numbers and some relation between these parameters in that family of hypergraphs.
Our work generalizes previous results concerning generalized power hypergraphs. 

\end{abstract}

\section{Introduction}
\hspace{0.5 cm}
A {\it (finite)  hypergraph} is a pair $H=(V(H), E(H))$ consisting of a (finite) non-empty set $V(H)$  and a collection $E(H)$ of non-empty subsets of $V(H)$. The elements of $V(H)$ are called {\it vertices} and the elements of $E(H)$ are called {\it hyperedges}, or simply {\it edges} of the hypergraph. The cardinal of  $V(H)$  is named the {\it order} of $H$. A {\it $k$-uniform } hypergraph is a hypergraph such that each edge consists of $k$ vertices.  The {\it rank} of $H$ is the maximum size of its hyperegdes (see \cite{Berge} for details).

Two vertices of $H$, $u$ and $v$,  are {\it adjacent} if there is an edge $e$ such that $u, v \in e$. The set of adjacent vertices in $H$ to a vertex $v$ is denoted as $N_H (v)$, or as $N(v)$ for simplicity. The vertex $v$ is {\it incident} to the edge $e$ if $v\in e$. The {\it degree} of a vertex $v$ is the number of edges incident to  $v$.  All over the paper,  a  {\it graph} will be a 2-uniform hypergraph. 

For a simple graph $G$,  the notion of \textit{Berge-$G$} hypergraph can be found in \cite{Gerbner}. We consider the accurate definition from \cite{Lu}.
  Given a simple graph $G=(V(G),E(G)),$ a hypergraph $H=(V(H),E(H))$  is called a {\it Berge-G} if there are an injection $i: V(G) \mapsto  V(H)$ and  a bijection $f : E(G) \mapsto E(H)$  such that for each $e=uv \in  E(G)$  we have $\{i(u), i(v)\} \subseteq  f(e)$. In other words, a hypergraph $H$ is a Berge-$G$ if we can embed a distinct edge of $G$ into each hyperedge of $H$ to create a copy of the graph $G$ on $H$. From now on and for the sake of simplicity, we identify $i(v)$ with $v$ for all $v \in V(G)$, unless there is ambiguity.

Trivially, $E(G)$ and $E(H)$ have the same cardinality. Let us remark that, given a graph $G$,  a Berge-$G$ can be constructed by replacing each edge of $G$ with a hyperedge that contains it (been allowed to introduce new vertices). $G$ is called {\it the support graph of $H$} and vertices of $V(G)\cap V(H)$ are named {\it support vertices of $H$}. Let us observe that an infinite family of  Berge-$G$ hypergraphs is obtained from the graph $G$; this family is denoted ${\mathcal B}(G)$. 


In the literature, Berge hypergraphs predominantly appear in Ramsey-type problems \cite{Ramsey problems} and Turán-type problems \cite{Gerbner}. Current studies on hypergraphs also explore parameters derived from graph theory, such as the clique number (i.e., the order of the largest complete subgraph), diameter, chromatic number, independence number, and others; see, e.g., \cite{Portugal}.

In spite of the concept of a hypergraph is relatively recent \cite{Berge}, it constitutes a versatile tool, due to its capacity to model complex relationships. In fact, hypergraphs are already being applied in various fields since many years ago, such as data mining \cite{Gunopulos}, artificial intelligence \cite{Eiter} and game theory \cite{Kumabe}, among others. However, the application of these concepts still offers significant room for further development.

Next, let us recall the parameters of our study. 

A {\it matching} in  a hypergraph $H$ is a set of disjoint hyperedges. The {\it matching number} of $H$, $\nu (H)$, is the maximum size of a matching in $H$. 

Given a hypergraph  $ H=(V, E),$ a subset $T\subset V$ is a {\it transversal} (or a {\it vertex cover}) of $H$ if $T$ has nonempty intersection with every hyperedge of $H$. The {\it transversal number} of $H$, $\tau(H),$ is the minimum size of a transversal of $H$. Transversals in hypergraphs are extensively studied in the literature (see, for example, \cite{Alon, Henning-Yeo}).

Given a hypergraph  $ H=(V, E),$   $D \subset V$  is   a {\it dominating set} of $ H$ if for every $v\in V - D$ there exists $u\in D$ such  that $u$ and $v$ are adjacent. The minimum cardinality of a dominating set of $H$, $\gamma(H)$, is  its {\it dominating number}. 
As it is indicated in the \textit{Introduction} of \cite{Shan}, dominating sets are important objects in communication networks, as they represent parts from which the entire network
can be reached directly. Messages can then be transmitted from sources to destinations via such a ‘‘backbone’’ with suitably
chosen links. To the best of our knowledge, the concept of domination in hypergraphs was
introduced by Acharya \cite{Acharya} and has been surveyed profusely further in, e.g., \cite{Acharya 2, Henning 2012, Jose-Tuza}.

\section{Dilations}
\hspace{0.5 cm}As a particular case of Berge-$G$ hypergraph, we introduce the  notion of dilation as follows.

\begin{definition}

Let  $G = (V(G),E(G))$ be a graph and $k \geq  3$. For each pairwise
adjacent vertices $v_i,\, v_j \in V(G)$, let $s_i$ and $s_j$ be two positive integers related to $v_i$ and $v_j$, respectively, and such that $2 \leq s_i + s_j \leq k$. A {\rm dilation} of $G$ is a hypergraph $H = (V(H),E(H))$ of rank $k$ whose vertex set is 
$$V(H) = (\bigcup_{v_i \in V(G)} {\bf v_i}) \cup(\bigcup_{e \in E(G)}{\bf e})$$ and whose hyperedges set is $$E(H)=\{ {\bf v_i}\cup  {\bf v_j} \cup {\bf e} : e = \{v_i, v_j\} \in E(G)\},$$
where ${\bf v_t}$ is a $s_t$-set containing $v_t$ ($t\in\{i, j\}$) and ${\bf e}$ is a set of size not greater than   $k-s_i-s_j$ which corresponds to the edge $e$. A vertex in ${\bf v_t}-v_t$ is called {\rm copy vertex} and a vertex in ${\bf e}$ is called {\rm additional vertex.}

\end{definition}	

\begin{notation}
{\rm The family of all dilations of a graph $G$ is denoted $\Gamma(G)$. 
The set of dilations without additional vertices in their hyperedges is denoted  $\Gamma_0(G)$ and  $\Gamma_1(G)$ denotes the set of dilations with at least one additional vertex in every hyperedge.} 
\end{notation}
Clearly, for any graph $G$, we get $\Gamma_0(G)\cap \Gamma_1(G) =\emptyset$, and for $G\neq K_2$, it is verified $\Gamma_0(G)\cup \Gamma_1(G) \subsetneq  \Gamma(G).$

\begin{figure}[htb]
	\centering
 \includegraphics[width=0.65\textwidth]{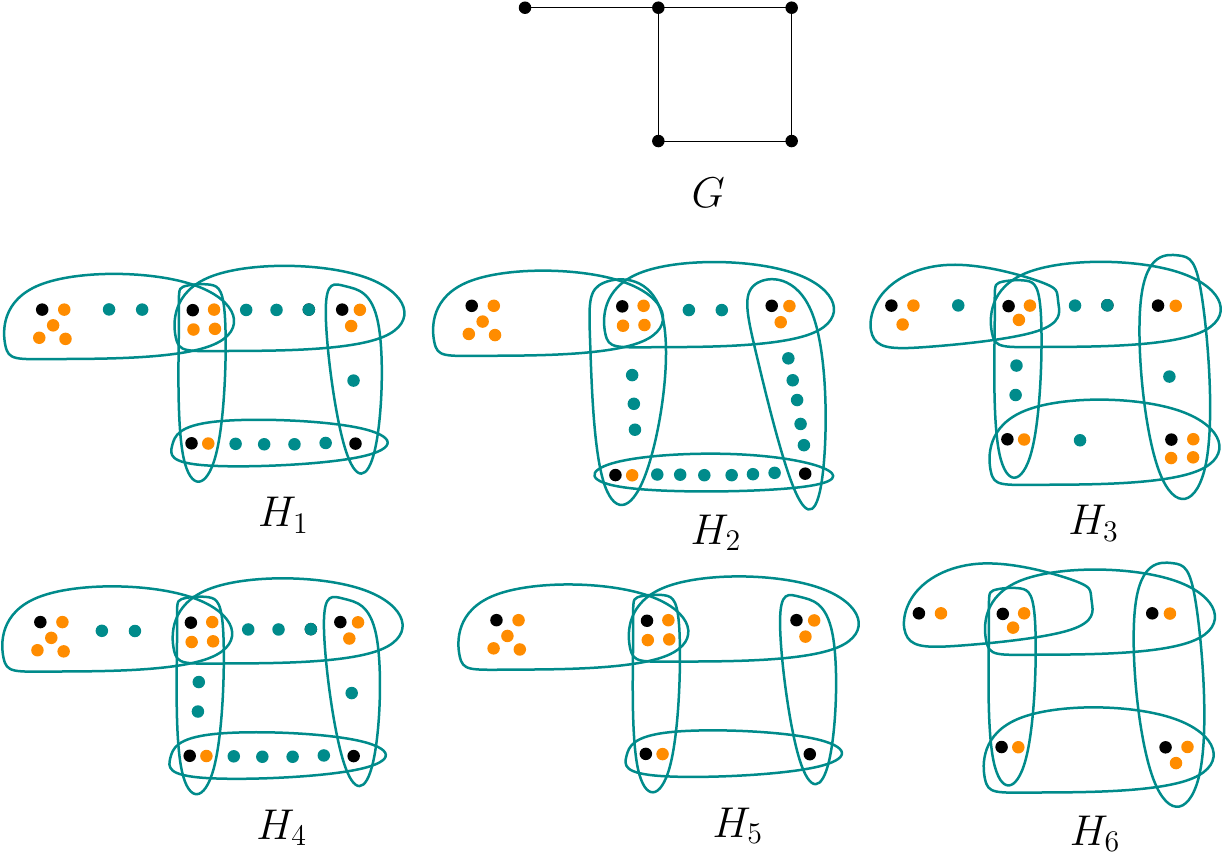}
	\caption{$H_1$, $H_2 \in \Gamma(G) \setminus (\Gamma_0(G) \cup \Gamma_1 (G))$, $H_3$, $H_4 \in \Gamma_1(G)$, $H_5$, $H_6 \in \Gamma_0(G)$.}\label{Dilataciones}
\end{figure}

As a particular case of a  Berge-$G$,  for any $k\geq 3$ and $1\leq s\leq \frac{k}{2}$, the {\it generalized power} of $G$, denoted by $G^{k,s}$, is defined in \cite{KhanFan2015} as the $k$-uniform hypergraph with   vertex set $V(G^{k,s})=(\bigcup_{v\in V} {\bf v})\cup (\bigcup_{e\in E} {\bf e})$, and   edge set $E(G^{k,s})=\{{\bf u}\cup{\bf v}\cup{\bf e}\,:\, e=\{u,v\}\in E \}$, where ${\bf v}$ is a $s$-set containing $v$, and ${\bf e}$ is a $(k-2s)$-set corresponding to $e$. Each vertex in  ${\bf v}-v$ is said to be {\it a copy of $v$.}  For $s=1$, $G^{k,1}=G^k$ is the {\it kth}-power hypergraph of $G$ and there is no copy of $v$ in $G^k$. This notion is defined as \textit{$k$-uniform expansion} of $G$ in \cite{Gerbner}.

It is readily checked that $G^{k,s}\in {\mathcal B}(G)$ for any $k\geq 3$ and $1\leq s\leq \frac{k}{2}$. Moreover,  for any $k\geq 3,$ it is also verified that:
\begin{itemize}
    \item If $1\leq s\leq \frac{k}{2}$ and $G\neq K_2$, then $G^{k,s}\in  \Gamma(G)\subsetneq {\mathcal B}(G)$ and $G^{k,\frac{k}{2}}\in  \Gamma_0(G)$.
    \item If $1\leq s< \frac{k}{2}$, then    $G^{k,s}\in  \Gamma_1(G).$
\end{itemize}  

This way, the concept of  dilation generalizes the power hypergraph ones. Notice that generalized power hypergraphs are uniform dilations, but they are not the only ones: see $H_3$ and $H_6$ in Figure \ref{Dilataciones}.\\ 

The basic properties of dilations are listed below.

\begin{proposition}\label{basicprop}
Let $G$ be a graph and $H \in \Gamma(G)$. Then:
\begin{enumerate}[label=\alph*)] 
    \item Each hyperedge of $H$ contains exactly two support vertices of $G$. 
    \item Two support vertices are adjacent in $H$ if, and only if,  they are adjacent in $G$.
    \item Two hyperedges of $H$ are disjoint if, and only if, the corresponding edges in $G$ are disjoint.
    \item $H$ is connected if, and only if, $G$ is connected.
\end{enumerate}
\end{proposition}

\begin{figure}[htb]
	\centering
 \includegraphics[width=0.5\textwidth]{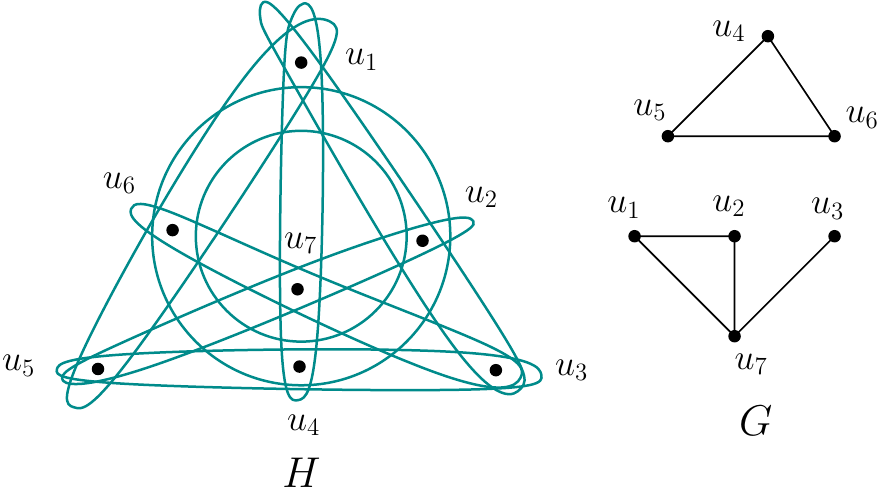}
	\caption{Fano plane \cite{Henning-Yeo} is a Berge-$G$.}\label{Plano de Fano}
\end{figure}

From the definition, any Berge-$G$ hypergraph  contains at least two support vertices in each hyperedge. Furthermore, for  Berge-$G$ hypergraphs, there are examples showing that the properties \textit{b)}, \textit{c)} and \textit{d)} are sufficient conditions, but not necessary ones, as we see, respectively, in Figures \ref{Contraejemplo adyacencia}, \ref{Contraejemplo hiperaristas no disjuntas} and \ref{Plano de Fano}. In particular, each support vertex of a Berge-$G$ $H$ has a degree in $H$ not lower than its degree in $G$; in dilations, due to the property \textit{b)}, degrees of support vertices are preserved.

\begin{center}
\begin{minipage}{0.48\textwidth}
    \centering
    \includegraphics[width=\textwidth]{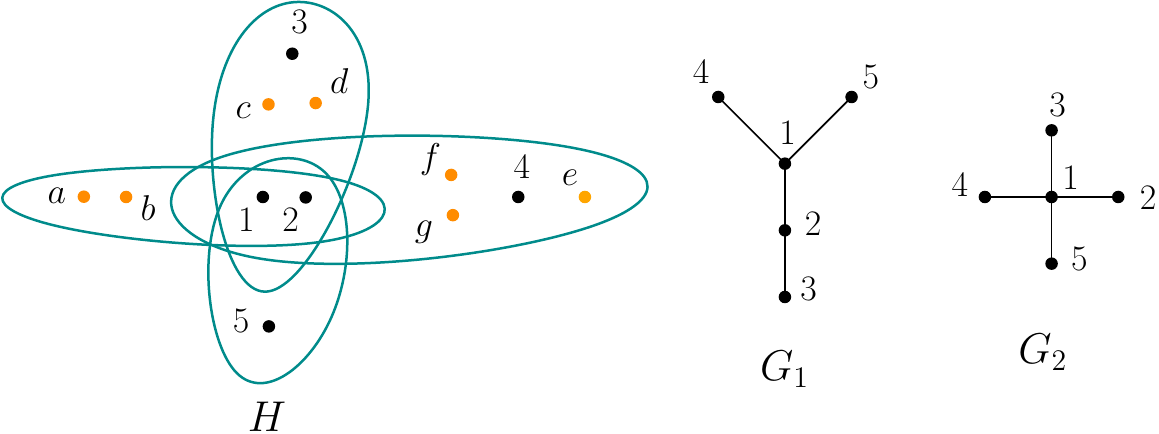}
    \captionof{figure}{$H \in \mathcal{B}(G_1) \cap \mathcal{B}(G_2)$.} 
    \label{Contraejemplo adyacencia}
\end{minipage}\hfill
\begin{minipage}{0.48\textwidth}
    \centering
    \includegraphics[width=\textwidth]{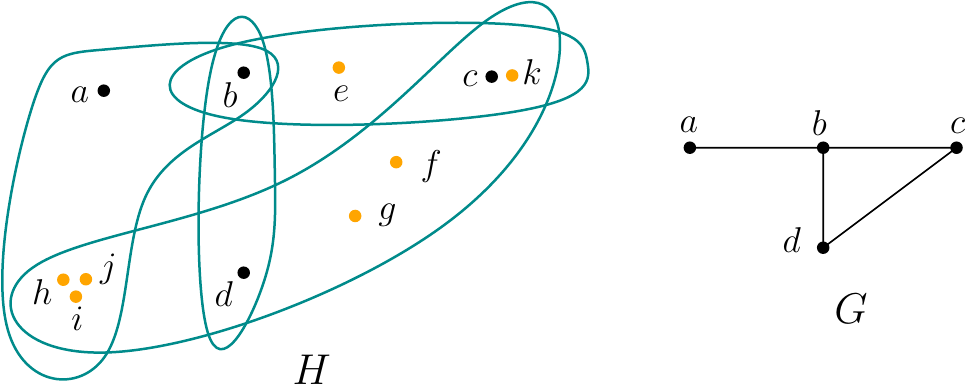}
    \captionof{figure}{$H \in \mathcal{B}(G)$.}
    \label{Contraejemplo hiperaristas no disjuntas}
\end{minipage}
\end{center}

 \section{Hereditary properties}
 
\hspace{0.5 cm}
For any graph $G$ with no vertex of degree $0$, it is known that 
$\gamma(G)\leq \nu(G)\leq \tau(G)$. For any hypergraph $H$, it is clear that $\gamma(H)\leq \tau(H)$ and $\nu(H)\leq \tau(H)$ hold. However, the inequality $\gamma(H)\leq \nu(H)$ does not hold for the general case.

To find out hypergraphs where  the corresponding parameters equalities are reached in the whole set of hypergraphs is a NP-hard problem \cite{Arumugan}.  This fact leads us to look for characterizations in a more restrictive family of hypergraphs.

As a starting point, we study properties of domination, matching and transversal numbers of $G$ that are inherited by the Berge-$G$ hypergraphs.  We get the following.

 \begin{theorem}\label{prop:igualdades}
 	
 	Let $G=(V,E)$ be a simple graph. Then:
 	\begin{enumerate}[label=\alph*)]
 	\item $\nu(H)\leq \nu(G)$ and $\tau(H)\leq \tau(G)$, for all $H \in \mathcal{B}(G)$.
        \item $\nu(H)=\nu(G)$, $\tau(H)=\tau(G)$ and $\gamma(G) \leq \gamma(H) \leq \tau(G),$ for all $H \in \Gamma(G)$.
 		\item $\gamma(H)=\gamma(G),$ for all $H \in \Gamma_0(G)$. 
 		\item $\gamma(H)=\tau(G),$ for all $H \in \Gamma_1(G)$. 
 	\end{enumerate}
 \end{theorem}
\begin{proof}
    
     Statement {\it a)}   is due to definitions.
    
     $\nu(H)=\nu(G)$ for all $H \in \Gamma(G)$ because of {\it a)}  and Proposition \ref{basicprop} \textit{c)}. 

     For {\it b)} , it is enough to prove $\tau(H)\geq \tau(G)$ for all $H \in \Gamma(G)$.
     Let $\mathcal{T}$ be a minimal transversal of $H$.
Each copy vertex of $\mathcal{T}$ can be exchanged by the support vertex of which it is a copy; suppose, then, that $\mathcal{T}$ contains no copy vertices. If $\mathcal{T}$ contains no additional vertices, then $\mathcal{T}$ is a subset of vertices of $G$, which implies that $\mathcal{T}$ is a transversal of $G$. Otherwise, it can be found another minimal transversal of $H$ with no additional vertices, $\mathcal{T'}$, by replacing each of these by one of its two adjacent support vertices.

Similarly, it can be proven that $\gamma(H)\geq \gamma(G)$ for all $H \in \Gamma(G)$. The other inequality derives from  $\gamma(H)\leq \tau(H)$ and, for all $H \in \Gamma(G)$, we have proved $\tau(H)=\tau(G)$.

Let us see {\it c)}. From {\it b)}, it is enough to prove $\gamma(G)\geq \gamma(H)$ for all $H \in \Gamma_0(G)$, what it is immediate since any minimum dominating set in $G$ is a dominating set in $H$.

Finally, let us prove {\it d)}. It is clear that $\gamma(H) \leq \tau(H)$ for any hypergraph $H$. In order to check $\gamma(H) \geq \tau(H)$, we claim that any domination set of $H$ is a transversal of $H$. Suppose not, then let $D$ be a dominating set in $H$ which is not a transversal of $H$. Therefore, there exits a hyperedge $e$ with no vertex in $D$. But every hyperedge in $H$ contains at least one additional vertex (whose degree is 1), then $H$ has at least one vertex (of $e$) which is not dominated by $D$, contradicting  the property of being $D$ a dominating set. Thus, we get $\gamma(H) \geq \tau(H)$.
\end{proof}

The equalities of statements {\it b)}, {\it c)} and {\it d)} will be referred as \textit{hereditary properties}.\\

To end this section, we add that it is not difficult to find examples of dilations in $\Gamma(G) \setminus (\Gamma_0(G) \cup \Gamma_1 (G))$, through the support graph $G$, whose domination number is $\gamma(G)$ or $\tau (G)$.\\



 \section{Extremal dilations}

 \hspace{0.5 cm}
 We extend some results related to connected generalized power hypergraphs proved in \cite{Dong} and present them within the context and terminology of dilations. In fact, this improving is quite direct since proofs also work when the condition of uniformity on generalized  power hypergraphs is replaced by the conditon of {\it  each hyperedge contains at least one additional vertex}. Therefore, we can state the following.




 \begin{theorem} \label{teoremaco}
 	Let $G$ be a graph. Then:
 	\begin{enumerate}[label=\alph*)]
 		\item $\nu(H)\leq \gamma(H) \leq 2 \nu(H),$ for all $H \in \Gamma_1(G)$.
 		\item $\gamma(H) \leq \nu(H),$ for all $H \in \Gamma_0(G)$.
 	\end{enumerate}
 \end{theorem}
 \begin{proof}
     If \( G \) is \( K_2 \), then the inequalities hold trivially, since  
\( \gamma(H) = \nu(H) = 1 \) for every dilation \( H \in \Gamma(K_2) \). From now on, let us assume that \( G \)  
is not \( K_2 \).  

Let \( H \in \Gamma_1(G) \). 
$\tau (G) \leq 2 \nu(G)$, as vertices of any maximum matching of $G$ form a transversal of $G$. So, by Theorem \ref{prop:igualdades}, it follows that \( \gamma(H) \leq 2\nu(H) \). Besides, from the general inequality $\nu(H) \leq \tau(H)$ and  Theorem \ref{prop:igualdades}, we get $\nu(H) \leq \gamma(H)$. 

On the other hand, if \( H \in \Gamma_0(G) \), it follows, from $\gamma(G) \leq \nu(G)$ and from Theorem \ref{prop:igualdades}, that \( \gamma(H) \leq \nu(H) \).
 \end{proof}
  
  In this context, the hypergraph $H$ is said to be {\it extremal for the domination number} if it verifies one of the  equalities $\gamma(H) = \nu(H)$ or $\gamma(H) = 2 \nu(H).$ The characterization of such extremal generalized power hypergraphs for the domination number is already introduced by Dong et al. in  \cite{Dong}; unfortunately, it is incomplete, as we explain in \cite{Arxiv}. To ease the reading, we include the details in Appendix \ref{appendix}. 

Studying such extremal hypergraphs by differentiating the  cases of $\Gamma_1(G)$ and $\Gamma_0(G)$ seems to be a good approach.
 

\subsection{Extremal dilations in $\Gamma_1(G)$}\label{subsección 1}



	\hspace{0.5 cm}
    König's Theorem \cite{Konig} affirms that any bipartite graph $G$ satisfies $\tau(G)=\nu(G).$ However, one can find non-bipartite graphs verifying this equality (see the graph in the Example \ref{ex:non bipartite}). A graph $G$ satisfying $\tau(G)=\nu(G)$ is named {\it König-Egerv\'ary graph} or said to {\it have the König- Egerv\'ary property} (KEG for short). 
    These graphs have been widely studied in the literature and several characterizations of them are given in terms of forbidden subgraphs (see \cite{Bonomo} and references there in).

    

	Hence,  knowing the family of KEG  will lead us to complete the correct characterization of power and generalized power  hypergraphs with equal domination and matching numbers. In addition, we extend that characterization for the family of dilations in $\Gamma_1(G)$. Also, we spread the characterization in \cite{Dong} of power and generalized power hypergraphs whose domination number is twice its matching number. Definitely, we get the following.
 
 \begin{theorem}\label{Mejorando con Gamma1}  
 	For  $H \in \Gamma_1(G)$ the following statement hold.
 	\begin{enumerate}[label=\alph*)]
 		\item $\gamma(H) = 2 \nu(H)$ if, and only if, $G=K_{2 \nu(H)+1}$.
 		\item 	$\gamma(H) = \nu(H)$ if,  and only if, $G$ is KEG.
 	\end{enumerate}
 \end{theorem}
 
 \begin{proof}
     From the computation of $\gamma$ and $\nu$ for the complete graphs with an odd number of vertices and from Theorem \ref{prop:igualdades}, we get the sufficient condition of {\it a)}.  
     Also, Gallai \cite{Gallai} proved that a graph $G$ verifies $\tau(G) = 2 \nu(G)$ if, and only if, $G=K_{2 \nu(G)+1}$; hence, by Theorem \ref{prop:igualdades}, we get the necessary condition in {\it a)}.
     

    For {\it b)}, let $H$ be a dilation in $\Gamma_1(G)$ such that $\gamma(H)=\nu(H)$. Then, by Theorem \ref{prop:igualdades}, $\tau(G)=\nu(G)$ holds. 
    Conversely, let us consider a hypergraph $H \in \Gamma_1(G)$ whose support graph is a KEG. From definition of KEG and from  Theorem \ref{prop:igualdades}, we get $\gamma(H)=\nu(H)$. 
 \end{proof}

	
\begin{example}\label{ex:non bipartite}
		{\rm Let $G = C_p \vee C_q$ denotes the   union of two cycles joined by precisely a common vertex, $p, q \geq 3$, being $p$ or $q$ an odd number. The following identities hold.

	\begin{itemize}
			\item $\nu(H) = \left \lfloor \frac{p}{2}  \right \rfloor + \left \lfloor \frac{q}{2}  \right \rfloor$. 
			
			\item $ \tau(H) = \left \lceil \frac{p}{2} \right \rceil + \left \lceil \frac{q}{2} \right \rceil - 1.$
	\end{itemize}	
		From Theorem \ref{prop:igualdades}, we get  $\gamma(G^{k,s})=\tau(G^{k,s})=\tau(G)$ for $1\leq s< k/2.$ Set $n=\nu(G)= \left \lfloor \frac{p}{2}  \right \rfloor + \left \lfloor \frac{q}{2}  \right \rfloor,  $ and let us distinguish two cases.
		
		\begin{enumerate}
			\item If $p$ and $q$ are  odd numbers, then $\left \lceil \frac{p}{2} \right \rceil + \left \lceil \frac{q}{2} \right \rceil - 1=n+1$ and, therefore
			$$\gamma(G^{k,s})=n+1=\nu(G)+1=\nu(G^{k,s})+1.$$
			\item If $p$ is even and $q$ is odd, then $\left \lceil \frac{p}{2} \right \rceil + \left \lceil \frac{q}{2} \right \rceil - 1= n$, and hence 
			$$\gamma(G^{k,s})=n= \nu(G^{k,s}).$$
		\end{enumerate} 
	
	Let us remark that every graph  of the form $C_p \vee C_q$, being $p$  even and $q$  odd, is not bipartite while  the  co\-rres\-ponding generalized power hypergraph is extremal for the equality of domination and matching numbers. 
	}
	\end{example}

\subsection{Extremal dilations in $\Gamma_0(G)$} \label{subsección 2}
\vspace{0.1 cm}

\hspace{0.5 cm}Let us begin with the characterization of connected graphs $G$   with no vertex of degree 1 that satisfy $\gamma(G)=\nu(G)$.

\begin{lemma} [\cite{RanderarthVolkmann1999}] \label{teoremadeladiscordia}  
	Let $G=(V(G),E(G))$ be a connected, bipartite graph with a minimum degree not less than 2 and such that $V(G) = V_1 \cup V_2$ in  such a way that $1 \leq \left| V_1\right| \leq \left| V_2\right|$. Then, $\gamma(G)=\nu(G)$ if, and only if,  
	for any two distinct vertices $x_1, x_2 \in V_1$ adjacent to a common vertex, 
	there are at least two distinct vertices $y_1, y_2\in V_2$ such that $y_i$ is adjacent uniquely to the vertices $x_1$ and $x_2$, for each $i \in \{1, 2\}$. Moreover, for these graphs,  $\gamma(G)=\nu(G)= |V_1|$  holds.  
\end{lemma}  

Let us denote by $\mathcal{G}_{\geq 2}^B$ the family of graphs described in Lemma \ref{teoremadeladiscordia}. 

\begin{figure}[htb] 
	\begin{center}  
		\includegraphics[width=4cm]{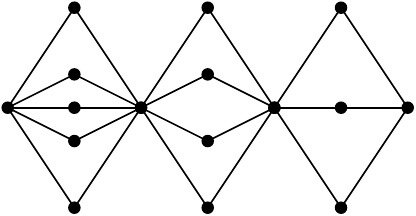}  
	\end{center}  
	\caption{Graph belonging to $\mathcal{G}_{\geq 2}^B$.} \label{Imagen_Grafo_G_2^B}  
\end{figure}

\begin{figure}[htb] 
	\begin{center}  
		\includegraphics[width=12.7cm]{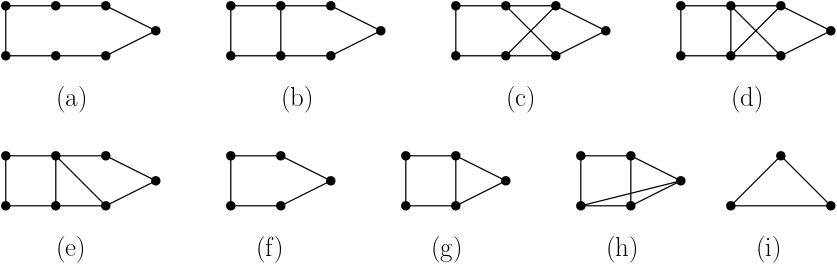}  
	\end{center}  
	\caption{The graphs in $\mathcal{G}_{\geq 2}^{NB}$ \cite{Kano}.} \label{Imagen_9_fantásticos}  
\end{figure}

\begin{lemma} [\cite{RanderarthVolkmann1999}] \label{teoremadeladiscordia2}  
	Let $G$ be a connected, non-bipartite graph of order $n$ with 
    minimum degree not less than 2. Then, $\gamma(G)=\nu(G)$ if, and only if, $G \in \mathcal{G}_{\geq 2}^{NB}$, the family consisting of the nine graphs depicted in Figure \ref{Imagen_9_fantásticos}.  
	Moreover, $\gamma(G)=\nu(G)= \left\lfloor \frac{n}{2} \right\rfloor$ if $G \in \mathcal{G}_{\geq 2}^{NB}$.  
\end{lemma}  

For the case of connected graphs $G$ with a minimum degree equals to 1, we will describe the family $\mathcal{G}_1$ defined by Kano et al. in \cite{Kano}.

To make the descriptions easier to follow, some graph definitions are introduced.  
  
\textit{End(G)} denotes the set of leaves (1-degree vertices) of $G$; an edge incident to a leaf is called a pendant edge. A vertex adjacent to a leaf is called a stem, and $Stem(G)$ denotes the set of stems of $G$. The corona $G \circ K_1$ of a graph $G$ is the graph obtained from $G$ by adding a pendant edge to each vertex of $G$. A connected graph $G$ of order at least 3 is a generalized corona if $V(G) = End(G) \cup Stem(G)$.  

\begin{figure}[htb]  
	\begin{center}  
		\includegraphics[width=7cm]{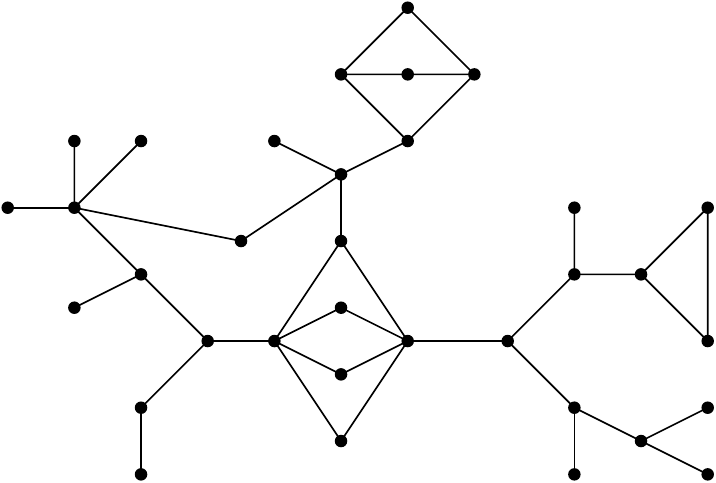}  
	\end{center}  
\caption{Example of a graph in $\mathcal{G}_1$ that is neither $K_2$ nor a generalized corona.} \label{Imagen_Grafo_de_G_1}  
\end{figure}

The family $\mathcal{G}_1$ (see Figure \ref{Imagen_Grafo_de_G_1} as an example) consists of $K_2$, generalized coronas and graphs whose connected components $G_j$ of $G - (End(G) \cup Stem(G))$, for $j \geq 1$, satisfy one of the following conditions:  
\begin{enumerate}   
	\item [i)] $G_j$ is the trivial graph.
	\item [ii)] $G_j$ is a connected bipartite graph with partition $V_1$ and $V_2$ where $1 \leq \left| V_1\right| < \left| V_2\right|$. Let $U_{G_j} = V(G_j)\cap N_G(Stem(G))$. Then, $\varnothing \neq U_{G_j} \subseteq V_2$ and for any two distinct vertices $x_1, x_2 \in V_1$ adjacent to a common vertex of $V_2$, there exist at least two distinct vertices $y_1,y_2 \in V_2 \setminus U_{G_j}$ such that $N_{G_j}(y_i)=\left\{x_1,x_2 \right\}$ for each $ i \in \left\{1,2 \right\}$.
	\item [iii)] $G_j$ is isomorphic to (f), (g), (h), or (i) from Figure \ref{Imagen_9_fantásticos}, and $\gamma(G_j - V_1)=\gamma(G_j)$ for all $ \varnothing \neq V_1 \subseteq U_{G_j} \subset V(G_j)$, where $U_{G_j} = V(G_j)\cap N_G(Stem(G))$.     
\end{enumerate}

It must be pointed out that we use a new unified notation for the families $\mathcal{G}_{\geq 2}^{B}$, $\mathcal{G}_{\geq 2}^{NB}$ and  $\mathcal{G}_1$.

The family $\mathcal{G}_1$ characterizes graphs with minimum degree 1 where the domination and matching numbers are equal.  
\begin{lemma} [\cite{Kano}] \label{Kanoresuelve}  
	A connected graph $G$ with a minimum degree equals 1 satisfies $\gamma(G)=\nu(G)$ if, and only if, $G \in \mathcal{G}_1$.  
\end{lemma}  

As a consequence of Lemmas \ref{teoremadeladiscordia}, \ref{teoremadeladiscordia2} \ref{Kanoresuelve} and Theorem \ref{prop:igualdades}, we get a characterization of extremal dilations in $\Gamma_0(G)$.

\begin{theorem}\label{Mejorando con Gamma0}  
	A connected dilation $H \in \Gamma_0 (G)$ satisfies $\gamma(H) = \nu(H)$ if, and only if, $G \in \mathcal{G}_{\geq 2}^B \cup \mathcal{G}_{\geq 2}^{NB} \cup \mathcal{G}_1$.  
\end{theorem}  

As a conclusion of this section, we want to emphasize that Theorems \ref{Mejorando con Gamma1} and \ref{Mejorando con Gamma0} not only improve the characterization given by Dong et al.(\cite{Dong}, Thm.~3.4), but they spread  it to a wide family of hypergraphs, namely $\Gamma_1(G)$ and $\Gamma_0(G)$. 

 \section{Non-extremal dilations}
 \hspace{0.5 cm}
 Now we present two results ensuring the existence of dilations for each non-extremal value of $\gamma(H)$ for $H \in \Gamma_0(G) \cup \Gamma_1(G)$; that is $\nu(H)< \gamma(H)< 2 \nu(H),$ for all $H \in \Gamma_1(G)$;
 		or else $\gamma(H) < \nu(H),$ for all $H \in \Gamma_0(G)$. To start with, we present families of graphs for the  cases $\gamma(H) = \nu(H) +1$,  $\gamma(H) = 2\nu(H) -1,$ and $\nu(H) = n$ but $\gamma(H) = 1.$

 \begin{lemma}\label{lemaNoExtremales}
     For any integer $n \geq 2$, there is a graph $G$ such that:
    \begin{enumerate}[label=\alph*)]
 		\item $\nu(H) = n$ and $\gamma(H) = \nu(H) +1 = n+1$ for any $H \in \Gamma_1(G)$.
 		\item $\nu(H) = n$ and $\gamma(H) = 2\nu(H) -1 = 2n-1$ for any $H \in \Gamma_1(G)$.
        \item $\nu(H) = n$ and $\gamma(H) = 1$ for any $H \in \Gamma_0(G)$.
    \end{enumerate}
 \end{lemma}
 \begin{proof}
     For {\it a)}, let us take  the cycle graph $G=C_{2n+1}$, $n \geq 2$, and apply  Theorem \ref{prop:igualdades}.
     For {\it b)} and {\it c)}, we consider $G$ as the complete graph $K_{2n}$, $n \geq 2$, and apply Theorem \ref{prop:igualdades}.
 \end{proof}
 
 \begin{theorem}\label{th:nonextremal_2n}
 	Let  $n \geq 2$ be an integer number. For any $ m \in \{{n+1}, \dots,\\ {2n-1}\}$ there is a graph $G$ such that $\nu(H)=n$ and $\gamma(H)=m$, for all $H \in \Gamma_1(G)$. 
 \end{theorem}
\begin{proof}
According to Lemma \ref{lemaNoExtremales}, there are support graphs \( G \)  
for which \( \nu(H) = n \), with \( n \geq 2 \), and \( \gamma(H) = m \), where  
\( m \in \{n+1, 2n-1\} \), for every \( H \in \Gamma_1(G) \).  
Therefore, it suffices to construct dilations \( H \), by means of their support graph \( G \),  
for which \( \nu(H) = n \), with \( n \geq 2 \), and \( \gamma(H) \) takes integer values in the range  
\( [n+2, 2n-2] \).

It is easily checked that for any the complete graph \( G = K_{2n} = (V(G),E(G)) \) and  \( H \in \Gamma_1(G) \),  \( \nu(H) = n \). Recall that any set of $r$ vertices of $V(G)$, with \( 2 \leq r \leq n - 1 \), induces a complete subgraph $K_r \subset G$.

Let us also consider the graph \( G'_r = (V(G),E(G) \setminus E(K_r)) = K_{2n} \setminus E(K_r)\).

Let \( H'_r \in \Gamma_1(G'_r) \).  

There exists a labeling of \( V(G) \) such that the edges of the clique \( K_r \) do not belong to  
the cycle \( \{v_1, v_2, \dots, v_{2n}, v_1\} \). In this cycle of \( G'_r \), we can find a perfect matching,  
so, by Theorem \ref{prop:igualdades}, it follows that \( \nu(H'_r) = \nu(G'_r) = n \).  

The set of vertices of \( G \) that do not belong to the subgraph \( K_r \) forms a minimum transversal of \( G'_r \). Hence, by Theorem \ref{prop:igualdades}, \( \gamma(H'_r) = 2n - r \), which completes the proof.

 \end{proof}
 
 \begin{theorem}\label{th:nonextremal_n}
 	Let $n \geq 2$ be an integer number.  For any $m \in \{1, \dots, n-1\}$ there is a graph $G$ such that $\nu(H)=n$ y $\gamma(H)=m$, for all $H \in \Gamma_0(G)$.
 \end{theorem}
 \begin{proof}
According to Lemma \ref{lemaNoExtremales}, there exist support graphs \( G \) for which  
\( \gamma(H) = 1 \) and \( \nu(H) = n \), with \( n \geq 2 \), for every \( H \in \Gamma_0(G) \).  
By Theorem \ref{prop:igualdades}, it suffices to construct graphs whose matching number is  
\( n > 2 \) and whose domination number is an integer \( m \) in the range  
\( [2, n - 1] \).  

Let \( T_1  = C_3 \vee K_2\) and \( T_2 = K_2 \vee_{v_1} C_3 \vee_{v_2} K_2 \) be the graphs depicted in Figure \ref{Fano_}.  

\begin{figure}[htb]
	\centering
 \includegraphics[width=0.35\textwidth]{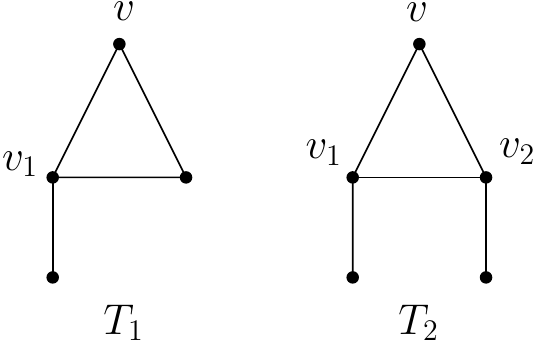}
	\caption{$T_1$ and $T_2$ graphs.}\label{Fano_}
\end{figure}

For every integer  $n>2$, let be $r \in \{1, \dots, \left \lfloor \frac{n-1}{2} \right \rfloor \}$.\\

Let us construct the graph $G(n,r)$ consisting of $n-2r$ copies of the graph $K_3$ and $r$ copies of the graph $T_2$ joined in a shared 2-degree vertex, $v$:
\begin{center}
$G(n,r) = C_3 \vee_v$ 
$\begin{matrix}
	n-2r) \\ 
	\ldots
\end{matrix}$
$\vee_v C_3 \vee_v T_2 \vee_v$ 
$\begin{matrix}
	r) \\ 
	\ldots
\end{matrix}$
$\vee_v T_2$ 
\end{center}

Let us also construct the graph $\widehat{G}(n,r)$ which is the union of  $n-2r$ copies of the graph $K_3$, $r-1$ copies of the graph $T_2$ and one graph $T_1$ joined in a shared 2-degree vertex, $v$: 

\begin{center}
$\widehat{G}(n,r) = C_3 \vee_v$ 
$\begin{matrix}
	n-2r) \\ 
	\ldots
\end{matrix}$
$\vee_v C_3 \vee_v T_2 \vee_v$ 
$\begin{matrix}
	r-1) \\ 
	\ldots
\end{matrix}$
$\vee_v T_2 \vee_v T_1$
\end{center}

\begin{figure}[htb]
	\centering
 \includegraphics[width=0.8\textwidth]{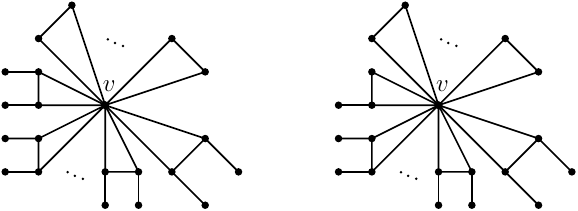}
	\caption{$G(n,r)$ and $\widehat{G}(n,r)$ graphs.}\label{Fano}
\end{figure} 

Firstly, since \( \nu(C_3) = 1 \) and \( \nu(T_2) = \nu(T_1) = 2 \), it follows that the matching number is  
\( \nu(G(n, r)) = 2r + (n - 2r) = n \). Similarly, we verify that  
\( \nu(\widehat{G}(n, r)) = 2r + (n - 2r) = n \).  

Secondly, due to the very construction of \( G(n, r) \), we have  
\( \gamma(G(n, r)) = 2r + 1 \), since the set of vertices consisting of \( v \) and the leaves  
of \( G(n, r) \) forms a minimum dominating set of \( G(n, r) \).  
Similarly, \( \gamma(\widehat{G}(n, r)) = 2r + 1 - 1 = 2r \), since \( \widehat{G}(n, r) \) has  
one fewer leaf than \( G(n, r) \).  

In summary, we have  
\( \nu(G(n, r)) = \nu(\widehat{G}(n, r)) = n \),  
\( \gamma(G(n, r)) = 2r + 1 \) and \( \gamma(\widehat{G}(n, r)) = 2r \). 

Therefore, for $n >2 $ and $m$ an odd integer in the range \( [2, n - 1] \), by  considering $r=\frac{m-1}{2}$ the graph $G (n, \frac{m-1}{2})$ verifies \( \gamma(G(n, r)) = m \). For $n >2 $ and $m$ an even integer in the range \( [2, n - 1] \), it suffices to consider $r=\frac{m}{2}$ and the graph $\widehat{G} (n, \frac{m}{2})$ verifies \( \gamma(\widehat{G}(n, r)) = m \), concluding the proof.
 \end{proof}

As a future work, we propose to characterize, by means of the support graph $G$, the families of dilations for which  each non-extremal value of $\gamma(H)$ given by Theorems \ref{th:nonextremal_2n} and \ref{th:nonextremal_n} is reached. In other words, to characterize graphs $G$ such that, respectively, its transversal number is $\nu(G) +r$ or its domination number is $\nu(G)-r$, for $r \in \{1, \dots, \nu(G)-1 \}$.

\appendix 
\section{Appendix: A comment to the characterization of extremal generalized power hypergraphs }\label{appendix}

\hspace{0.5 cm}
Next, we collect two principal results stated by Dong et al. in \cite{Dong}, in order to point out what seems to be the main confusion in the failed
characterization given in the following theorem.

\begin{theorem}\label{pseudoteoremas}{\rm(\cite{Dong}, Thm.~2.2 and Thm.~3.4)}
		Let $k\geq 3$. Any connected  generalized power hypergraph $G^{k,s}$ satisfies $\gamma(G^{k,s})= \nu (G^{k,s})$ if, and only if, $G$ is a bipartite connected graph, for $1 \leq s < \frac{k}{2}$, or if, for $s= \frac{k}{2}$, $G$ is a connected graph of the family ${\mathcal G}_1\cup{\mathcal G}_{\geq 2}^{NB}$.
\end{theorem}



	The proofs of both Theorems  2.2 and   3.4 in \cite{Dong} are quite similar for the case $1\leq s< k/2$, hence we  focus on reasoning only Theorem  2.2. For the necessary condition it is  affirmed that the equality of domination and matching numbers is raised only by bipartite graphs; however, it is also raised for other graphs, as we have seen in Example \ref{ex:non bipartite}. The misunderstanding is produced because they use   König's Theorem as a characterization of  graphs with equal transversal and matching numbers. But the fact is that the bipartiteness is only a sufficient condition, not a necessary one as we check from König's Theorem \cite{Konig} ($\tau(G)=\nu(G)$ if $G$ is a bipartite graph).

Theorem \ref{pseudoteoremas}, in the case $s = \frac{k}{2}$, establishes a sufficient condition for $\gamma(G^{k, \frac{k}{2}}) = \nu(G^{k, \frac{k}{2}})$, but it is not necessary.
The following example refutes it. 
		\begin{example}	
			{\rm Let us consider the  graph $G=K_{2,n}$, $n\geq 2$. Hence $\gamma(G^{k,\frac{k}{2}})=  2 = \nu(G^{k,\frac{k}{2}}).$ However, it is easy to check that $G$ does not belong to $\mathcal{G}_{\geq 2}^{NB} \cup \mathcal{G}_1$.}    
\end{example}

The analysis of \textit{Theorem} 3.4 in \cite{Dong} for the case  $s = \frac{k}{2}$ is tackled now (observe that the case $s = \frac{k}{2}$ does not appear in \textit{Theorem} 2.2 in \cite{Dong}). The proof of the necessary condition is based in necessary conditions of \textit{Lemmas} 3.2 and 3.3 in \cite{Dong}. However, \textit{Lemma} 3.2 in \cite{Dong}, as stated, is quite misleading.
We have included in subsection \ref{subsección 2} the complete characterization of graphs with equal matching and domination numbers summarized from \cite{Dong, Kano, RanderarthVolkmann1999}.

As a consequence of Lemmas \ref{teoremadeladiscordia}, \ref{teoremadeladiscordia2} and \ref{Kanoresuelve}, by  using Theorem \ref{prop:igualdades}, the characterization of generalized power hypergraphs $G^{k,\frac{k}{2}}$ with equal domination
 and matching numbers is immediately deduced. 
 Therefore, we provide an improvement by establishing a necessary and sufficient condition for the generalized power hypergraphs $G^{k, \frac{k}{2}}$, with even $k \geq 4$. Moreover, this characterization has been placed in Theorem \ref{Mejorando con Gamma0} in a broader context by being framed within the family of dilations $\Gamma_0(G)$.


\end{document}